\documentclass[10pt]{article}
\usepackage{calc}
\usepackage[all]{xy}
\usepackage[centertags]{amsmath}
\usepackage{latexsym}
\usepackage{amsfonts}
\usepackage{amssymb}
\usepackage{amsthm}
\usepackage{fancyhdr}
\usepackage{url}
\usepackage{times}
\usepackage{color}
\newtheorem{thm}{Theorem}[section]
\newtheorem{defn}[thm]{Definition}
\newtheorem{rem}[thm]{\bf Remark}
\newtheorem{rems}[thm]{\bf Remarks}
\newtheorem{prop}[thm]{Proposition}
\newtheorem{lem}[thm]{Lemma}

\usepackage{amsmath}
\allowdisplaybreaks
\usepackage[pdfpagelabels,plainpages=false]{hyperref}
\hypersetup{ linkcolor =blue,
 citecolor = cyan,
 urlcolor = blue,
colorlinks = true, bookmarksopen = true, bookmarksnumbered = true,
bookmarksopenlevel = {1}, linktocpage = true,
pdfstartview = FitH,
pdftitle = {Non-linear Case},
pdfsubject
= {article},
pdfauthor = {Habiba  KNANI},
pdfcreator = {pdftextify}, pdfproducer = {pdftextify}
}
\begin{document}

\title{Partial differential equations associated to non linear backward stochastic differential equations with Gaussian Volterra processes}
\author{H.Knani \thanks{%
Laboratoire de Math\'{e}matiques: D\'{e}terministe et Al\'{e}atoires, Hammam
Sousse, Tunisie et Université de Lorraine, CNRS, IECL, F-54000 Nancy, France}}

\maketitle

\begin{abstract}
 {\footnotesize  In this paper,  we generalize to Gaussian Volterra processes the existence and uniqueness of solutions for a class of non linear  backward stochastic differential equations (BSDE) and we establish the relation between the non linear BSDE  and the partial differential equation (PDE). A comparison theorem for the solution of the BSDE is proved and the continuity of its law is studied.
 }
  

\end{abstract}

\hrule
\maketitle

\section{Introduction}

\quad Backward stochastic differential equations (BSDE) driven by a brownian motion have been introduced by Bismut \cite{JMBism} in the linear case. Non linear BSDE have been studied first by Paradoux and Peng \cite{PaPe}. Since then BSDE have been of interest due to the connections with partial differential equations (PDE) and their applications, especially in mathematical finance, stochastic differential games and stochastic control.

\quad In this paper, we study the BSDE
\begin{equation}
\begin{cases}
 dY_t=-f(t,N_t,Y_t, Z_t)dt-Z_t\delta X_t, &\text{$t\in[t_0, T)$ with $t_0\geqslant 0$}\\
Y_T=g(N_T)\\
\end{cases}
\label{edsr}
\end{equation}
where $X=\{X_{t},0\leqslant t\leqslant T\}$ is a zero mean continuous Gaussian
process given by 
\begin{equation}
X_{t}=\displaystyle\int_{0}^{T}K(t,s)dW_{s},  
\label{X}
\end{equation}
where $W=\{W_{t},0\leqslant t\leqslant T\}$ is a standard Brownian motion
and $K:[0,T]^{2}\rightarrow 
\mathbb{R}
$ is a square integrable kernel, i.e. $\int_{[0,T]^{2}}K(t,s)^{2}dtds<+%
\infty .$ We assume that $K$ is of Volterra type, i.e, $K(t,s)=0$ whenever $%
t<s$. Usually, the representation (\ref{X}) is called a Volterra
representation of $X$. 
The kernel $K$ in (\ref{X}) defines a linear operator in $L^{2}([0,T])$
given by $(K\sigma )_{t}=\int_{0}^{t}K(t, s)\sigma _{s}ds,$ $\sigma \in
L^{2}([0,T]).$ 
 The process $(N_t, 0\leqslant t \leqslant T)$ is given by
\begin{equation}
N_t:=\int_{0}^{t}\sigma_s\delta X_s=\int_{0}^{t}(K_t^*\sigma)_s\delta W_s, \label{N}
\end{equation}
with $\sigma$ being a deterministic function and $K^*$ the adjoint operator of $K$ (\cite{AMN}, lemma 1) given by (\ref{Ketoile}) (see also \cite{DoKn}). $f$ is called the generator of the BSDE, $g(N_T)$ the terminal condition. \\
In \cite{AMN} $K$ is called \textit{regular }if it satisfies\textit{%
\ }

\textit{(H) } $\int_{0}^{T}\mid K\mid
((s,T],s)^{2}ds<\infty ,$ where $\mid K\mid ((s,T],s)$ denotes the total
variation of $K(.,s)$ on $(s,T]$.

We assume the following condition on $K(t,s)$ which is more restrictive than 
\textit{(H) }(\cite{AMN}, \cite{Coup.Maslk}):\\
\emph{(H1) \ }$K(t, s)$ is continuous for all $ 0<s \leqslant t<T $ and  continuously differentiable in the variable $t$ in 
${0<s< t<T}$,

\emph{(H2) \ }For some $0<\alpha, \beta <\frac{1}{2},$ there is a finite
constant $c>0$ such that
\begin{equation}
\begin{aligned} \Big| \frac{\partial K}{\partial t}(t, s)\Big| \leqslant
c(t-s)^{\alpha-1}\Big(\frac{t}{s}\Big)^{\beta},\,\text{for all $0<s<t<T.$}
\end{aligned}  \notag
\end{equation}
%

\quad Examples of Gaussian Volterra processes that
satisfy\ \emph{(H1) and (H2)} are multi-fractional Brownian motion (mbm), multi-fractional Ornstein-Uhlenbeck process and Liouville multi-fractional Brownian motion (\cite{DoKn}).\\
 \color{black}
\quad The covariance function of $X$ is given by 
\begin{equation}
R(t,s):=\mathbb{E}X_{t}X_{s}=\displaystyle\int_{0}^{inf(t,s)}K(t,u)K(s,u)du.
\label{Rh}
\end{equation}
\quad The aim of this paper is to study the nonlinear BSDE (\ref{edsr}) and we establish the relation to the associated partial differential equation (\ref{Edp}) that opens the possibility to solve the BSDE by means of classical or viscosity solutions of the PDE. This generalizes a resultat in \cite{HuPeng} obtained for fractional brownian motion. For the existence and  uniqueness for the solution of the BSDE (\ref{edsr}), essentially two methods were applied. The existence and uniqueness of the solution of (\ref{edsr}) is addressed in Theorem \ref{Thmsoledsr} by means of the associated PDE. Another proof will be treated in a separate paper without making reference to this PDE, but with probabilistic and functional theoretic methods. In this paper, we prove in Theorem \ref{connexionedp} that $Z_t=-\sigma_t\frac{\partial}{\partial x}u(t, x)$ for BSDE's with generators $f\in\mathcal{C}^{0, 1}([t_0, T]\times\mathbb{R}^3)$ of polynomial growth and under the injectivity hypothesis for the adjoint operator to $K.$ We discuss this hypothesis in Remark \ref{reminjectivite}. This hypothesis is satisfied for the mbf (with $H>\frac{1}{2}$) and comes from the preliminary Lemma \ref{lemmaab} that shows a kind of orthogonality between Lebesgue and divergence integrals. The proof of this Lemma generalizes the proof given by Y. Hu and S. Peng in 2009 for fractional brownian motion (\cite{HuPeng}). The proof of Lemma \ref{lemmaab} depends itself on Proposition \ref{finteg} where we show that, for any continuous function of exponential growth $h,$ $h(N_t)$ admits a representation as a divergence integral of the heat kernel operator evaluated for $N.$ We generalize in Theorem \ref{Comparisontheo} a comparison theorem, known for BSDE with respect to brownian motion (see for example \cite{ZhangI}), to the solution of the BSDE (\ref{edsr}) and study the continuity of the law of Y in Theorem \ref{continuityY}.

\quad Here is the organisation of the paper. In Section 2 we state the main results on the solution of (\ref{edsr}). In Section 3 we give some definitions and complements on the Skorohod integral with respect to Volterra processes.  Section 4 is devoted to the proofs of the results and contains other results of independent interest, like an Itô formula for $N$ proven in the framework of Malliavin calculus (Theorem \ref{Ito}) and a transfer formula (Proposition \ref{Transferformula}). 
\section{Statement of the main results}
We consider $X$ defined on the probability space $(\Omega, \mathcal{F}, P)$ and given by (\ref{X}). Let $\mathbb{F=\{}\mathcal{F}_{t}\subset \mathcal{F},$ $t\in \lbrack 0, T]\}$ the filtration generated by $X$ and augmented by the $P$-null sets. Let $N$ be given by (\ref{N}), where $\sigma$ is a bounded function on $[0,T]$ and suppose that $\frac{d}{dt}Var(N_t)>0$ for all $t\in(0, T).$


Let $t_{0}\geqslant 0$ be fixed, and denote by $\mathbb{L}^{2}(\mathbb{F},%
\mathbb{%
\mathbb{R}
})$ the set of $\mathbb{F}$-adapted $\mathbb{%
\mathbb{R}
}$-valued processes $Z$ such that $\mathbb{E}(\int_{t_{0}}^{T}\mid
Z_{t}\mid ^{2}dt)<\infty .$ We consider the non linear BDSE for the processes $%
Y=(Y_{t},$ $t\in \lbrack t_{0},T])\in \mathbb{L}^{2}(\mathbb{F},\mathbb{%
\mathbb{R}
})$ and $Z=(Z_{t},$ $t\in \lbrack t_{0},T])\in \mathbb{L}%
^{2}(\mathbb{F},\mathbb{%
\mathbb{R}
})$ given by

\begin{equation}
Y_{t}=g(N_{T})+\int_{t}^{T}f(s, N_s, Y_s, Z_s)ds+%
\int_{t}^{T}Z_{s}\delta X_{s},\text{ \ }t\in \lbrack t_{0},T],  \label{BSDE}
\end{equation}
We show that (\ref{BSDE}) is associated to the following
second order PDE with terminal condition%
\begin{equation}
\begin{cases}
\frac{\partial u}{\partial t} (t, x)=-\frac{1}{2}\frac{d}{dt}Var(N_t)\frac{\partial^2 u}{\partial x^2}(t,x)-f(t,x,u(t, x),-\sigma_t\frac{\partial u}{\partial x}(t, x)), \\
u(T,x)=g(x),\ \text{$(t,x)\in \lbrack t_{0},T]\times \mathbb{R}$}\\
\end{cases}
\label{Edp}
\end{equation}
The association of this PDE to the BSDE (\ref{BSDE}) is proven by means of the It\^o formula for the class of functions $F\in \mathcal{C}^{1, 2}([0, T]\times \mathbb{R})$ that satisfy, together with their partial derivatives, the growth condition\\
\begin{equation}
\begin{aligned} \max \Big( \Big | F(t,x)\Big |,\Big |\frac{\partial
F}{\partial t}(t,x)\Big |,\Big |{\frac{\partial F}{\partial x}}(t,x)\Big
|,\Big |{\frac{\partial^{2}F}{\partial x^{2}}(t,x)}\Big | \Big)
&\leqslant ce^{\lambda  x^{2}}, \end{aligned}  \label{f}
\end{equation}%
for all $t\in \lbrack 0,T]$ and $x\in \mathbb{R},$ where $c,\lambda $
are positive constants such that $\lambda<\frac{1}{4}(\underset{t\in \lbrack 0,T]}{\sup }Var(N_t))^{-1}$. This implies 
\begin{eqnarray}
\mathbb{E}\Big |F(t,N_{t})\Big |^{2} &\leqslant &c^{2}\mathbb{E}exp(2\lambda
|N_{t}|^{2})<\infty ,
\label{F2}
\end{eqnarray}
and the same property holds for $\partial /\partial tF(t,x),$ $\partial
/\partial xF(t,x)$ and $\partial ^{2}/\partial x^{2}F(t,x).$\\

The main results of this paper is stated below:
\begin{thm}
If $(\ref{Edp})$ has a classical solution $u(t, x)$ that satisfies (\ref{f}) with $F$ replaced by $u,$ then $(Y_t, Z_t):=(u(t, N_t), -\sigma_t\frac{\partial u}{\partial x}(t, N_t)),$ $t\in(t_0, T)$ satisfies (\ref{BSDE}) and $Y, Z\in\mathbb{L}^2(\mathbb{F},\mathbb{R}).$
\label{Thmsoledsr}
\end{thm}

\begin{rem}

\quad The mild (or evolution) solution of (\ref{Edp}) is given by
\begin{equation}
u(t, x)=\int G(T, t, x-y)g(y)dy-\int_{t}^{T}\int G(s, t, x-y)f(s, y, u(s, y),-\sigma
_{s}\frac{\partial }{\partial x}u(s, y))dyds,
\label{Gsolevol}
\end{equation}
where
$$
G(s, t, x)=(2\pi )^{-1/2}(Var(N_{s})-Var(N_{t}))^{-1/2}\exp \left( -\frac{x^{2}%
}{2(Var(N_{s})-Var(N_{t}))}\right).$$

\quad In fact, theorem 4.1 (\cite{DoKn}) applied to $n=1$
gives $u(t, x)=\int G(T, t, x-y)g(y)dy$ in the linear case and (\ref{Gsolevol})
follows by classical arguments for the mild form of a nonlinear PDE. 

\label{RemarkEDP}
\end{rem}
\begin{thm}
Suppose that $K_T^*$ is injective. Let $(\ref{BSDE})$ have a solution of the form $(Y_t=u(t, N_t), Z_t=v(t, N_t))$, where $u$ is as in Theorem \ref{Thmsoledsr} and $v \in \mathcal{C}^{1, 2}([0, T ] \times  \mathbb{R})$ satisfies (\ref{f}). 
Suppose that $f\in\mathcal{C}^{0, 1}([0, T ] \times  \mathbb{R}^3)$ and of polynomial growth. 
Then $v(t, x)=-\sigma_{t}\frac{\partial}{\partial x}u(t, x)$.
\label{connexionedp}
\end{thm}
\begin{rem}
If the PDE (\ref{Edp}) has a unique classical solution, then under the hypotheses of Theorem \ref{connexionedp} the BSDE (\ref{BSDE}) has a unique solution with $Y_t=u(t, N_t),$ $Z_t=v(t, N_t).$\\
In fact, let us suppose that $(\ref{BSDE})$ has two solutions  $(Y_t^{1}, Z_t^{1})$ and $(Y_t^{2}, Z_t^{2}):$ 
 $Y_t^{1}=u^1(t, N_t)$, $Z_t^{1}=-\sigma_t\frac{\partial}{\partial x}u^1(t, N_t)$ and $Y_t^{2}=u^2(t, N_t)$, $Z_t^{2}=-\sigma_t\frac{\partial}{\partial x}u^2(t, N_t)$. By uniqueness of solutions of (\ref{Edp}), $u^1(t, .)=u^2(t, .)$ and therefore, for $t\in[t_0, T]$ $Y_t^{1}= Y_t^{2}$, $Z_t^{1}=Z_t^{2}$.  
\end{rem}
\begin{thm} For $i=1, 2$ we consider the BSDE's 
\begin{equation}
Y_{t}^{i}=g^{i}(N_{T}^{i})+%
\int_{t}^{T}f^{i}(s, N_{s}^{i}, Y_{s}^{i}, Z_{s}^{i})ds+\int_{t}^{T}Z_{s}^{i}%
\delta X_{s}^{i},\text{ \ }0<t_{0}\leqslant t\leqslant T,\label{edsrmulti}
\end{equation}

where $N_{t}^{i}=\int_{0}^{t}\sigma _{s}^{i}\delta X_{s}^{i}$
and $c_{0}<\sigma ^{i}<C_{0}$ for some constants $c_{0}, C_{0}>0.$ Suppose that both BSDE's satisfy the hypotheses of Theorem \ref{Thmsoledsr}
and Theorem \ref{connexionedp}. If $f^{1}(\cdot , \cdot , y, z)\geqslant $ $
f^{2}(\cdot , \cdot , y, z)$ for all $y, z,$ and $g^{1}\geqslant g^{2},$ then $%
Y_{t}^{1}\geqslant Y_{t}^{2}$ $P-a.s.$ for all $t\in \lbrack t_{0}, T].$
\label{Comparisontheo}
\end{thm}
\begin{thm}
Suppose that $Y$ is a solution of the BSDE (\ref{BSDE}) and $u$ is a classical solution of the PDE (\ref{Edp}). If $\frac{\partial}{\partial x}u(t, N_t)\neq 0$ for all $t\in \lbrack t_{0}, T],$ the law of $Y$ is absolutely continuous.
\label{continuityY}
\end{thm}
\section{Preliminaries}
In this section, we recall important definitions and results concerning the Malliavin calculus for Volterra process. These results will be used to study the BSDE (\ref{BSDE}). \\

Let $\mathcal{E}$ be the set of step functions of $[0,T]$,
and let $K_{T}^{\ast }:\mathcal{E}\rightarrow L^{2}([0,T])$ be defined by 
\begin{equation}
(K_{T}^{\ast }\sigma )_{u}:=\int_{u}^{T}\sigma _s\frac{\partial K}{\partial
s}(s, u)ds.
\label{Ketoile}
\end{equation}%
\newline
%

\begin{rems} a) For $s>t,$ we have $(K_{T}^{\ast }\sigma
1_{[0,t]})_{s}=0,$ and we will denote $(K_{T}^{\ast }\sigma 1_{[0,t]})_{s}$
by $(K_{t}^{\ast }\sigma )_{s}$ where $K_{t}^{\ast }$ is the adjoint of the
operator $K$ in the interval $[0,t].$\newline
b) If $K(u, u)=0$ for all $u\in \lbrack 0,T],$ $(K_{T}^{\ast
}1_{[0,r]})_{u}=K(r, u)$ for $u<r.$ Indeed, if $u\leqslant r,$ we have\newline
\begin{equation*}
(K_{T}^{\ast }1_{[0,r]})_{u}=\displaystyle\int_{u}^{T}1_{[0,r]}(s)\frac{%
\partial K}{\partial s}(s, u)ds=\displaystyle\int_{u}^{r}\frac{\partial K}{%
\partial s}(s, u)ds.
\end{equation*}
\end{rems}
Therefore %
%
%
%
\begin{align*}
R(t,s)& =\mathbb{E}\Big[X_{t}X_{s}\Big]=\displaystyle%
\int_{0}^{inf(t,s)}(K_{T}^{\ast }1_{[0,t]})_{u}(K_{T}^{\ast }1_{[0,s]})_{u}du
\\
& =<K_{T}^{\ast }1_{[0,t]},K_{T}^{\ast }1_{[0,s]}>_{L^{2}([0,T])}.
\end{align*}
For $\sigma , \widetilde{\sigma }\in \mathcal{E}$ this may be extended to%
\begin{equation*}
X(\sigma ):=\displaystyle\int_{0}^{T}(K_{T}^{\ast }\sigma )_{s}dW_{s}\text{
\ and \ }E\Big[X(\sigma )X(\widetilde{\sigma })\Big]=<K_{T}^{\ast }\sigma
,K_{T}^{\ast }\widetilde{\sigma }>_{L^{2}([0,T])}.
\end{equation*}
\color{black}
\begin{defn}
Let $\mathcal{H}$ be the closure of the linear span of the indicator
functions $1_{[0,t]},$ $t\in \lbrack 0,T]$ with respect to the semi-inner product%
\begin{equation*}
<1_{[0,t]}, 1_{[0,s]}>_{\mathcal{H}}:=<K_{T}^{\ast }1_{[0,t]}, K_{T}^{\ast
}1_{[0,s]}>_{L^{2}([0,T])}.
\end{equation*}%
\label{defH}
\end{defn}
The operator $K_{T}^{\ast }$ is an isometry between $\mathcal{H}$ and a
closed subspace of $L^{2}([0,T]),$ and $\parallel \cdot \parallel _{\mathcal{%
H}}$ is a semi-norm on $\mathcal{H}.$ Furthermore, for $\varphi, \psi \in 
\mathcal{H}$, 
%
\begin{align*}
<K_{T}^{\ast }\varphi ,K_{T}^{\ast }\psi >_{L^{2}([0,T])}&=\int_{0}^{T}(K_T^*\varphi)_t(K_T^*\psi)_tdt\\
&=\int_{0}^{T}\int_{t}^{T}\varphi_r\frac{\partial K}{\partial r}(r, t)dr\int_{t}^{T}\psi_s\frac{\partial K}{\partial s}(s, t)dsdt\\
&=\int_{0}^{T}\displaystyle\int_{0}^{T}\Big(\displaystyle\int_{0}^{inf(r,s)}%
\frac{\partial K}{\partial r}(r, t)\frac{\partial K}{\partial s}(s, t)dt\Big)%
\varphi _{r}\psi _{s}dsdr.
\end{align*}

For further use let
\begin{equation}
\phi (r, s) :=\displaystyle\int_{0}^{inf(r, s)}\frac{\partial K}{\partial r}%
(r, t)\frac{\partial K}{\partial s}(s, t)dt,\text{ }r\neq s.
\label{phi}
\end{equation}
\begin{equation}
\widetilde{\phi }(r, s):=\displaystyle\int_{0}^{inf(r, s)}\Big | \frac{%
\partial K}{\partial r}(r, t)\Big | \Big |\frac{\partial K}{\partial s}(s, t)\Big |
dt,\text{ }r\neq s.
\label{phitilde}
\end{equation}
Note that $\phi (r, s)=\partial ^{2}/\partial s\partial r R(r, s)$ ($r\neq s)$ (%
$\phi $ may be infinite on the diagonal $r=s).$ 
%
%
%
%
%
%
%
Let $\mid \mathcal{H}\mid $ be the closure of the linear span of indicator
functions with respect to the semi-norm given by
\begin{align}
\parallel \varphi \parallel _{\mid \mathcal{H}\mid }^{2}&=\displaystyle%
\int_{0}^{T}\left( \int_{t}^{T}\mid \varphi _{r}\mid \Big | \frac{\partial K}{%
\partial r}(r, t)\Big | dr\right) ^{2}dt 
=2\displaystyle\int_{0}^{T}dr\displaystyle\int_{0}^{r}ds\widetilde{\phi }%
(r, s)\mid \varphi _{r}\mid \mid \varphi _{s}\mid.
\label{normephimodH}
\end{align}
We briefly recall some basic elements of the stochastic calculus of
variations with respect to $X$ given by (\ref{X}). We refer to \cite{IvPic} and \cite{Nualart} for a more
complete presentation. Let $\mathcal{S}$ be the set of random variables of
the form $F=f(X(\varphi _{1}),...., X(\varphi _{n}))$, where $n\geqslant 1$, $f : \mathbb{R}^n \rightarrow \mathbb{R}$ is a $\mathcal{C}^{\infty }$ -function such that $f$ and its partial derivatives have at most polynomial growth, and $\varphi _{1},..., \varphi _{n}\in \mathcal{H}$. The derivative
of $F$ 
\begin{equation*}
\begin{aligned} D^{X}F &:= \sum_{j=1}^{n}\frac{\partial f}{\partial
x_{j}}(X(\varphi_1),...., X(\varphi_n))\varphi_{j}. \end{aligned}
\end{equation*}%
\quad is an $\mathcal{H}$-valued random variable, and $D^{X}$ is a closable
operator from $L^{p}(\Omega )$ to $L^{p}(\Omega ;\mathcal{H})$ for all $%
p\geqslant 1$. We denote by $\mathbb{D}_{1,p}^{X}$ the closure of $\mathcal{S%
}$ with respect to the norm 
\begin{equation}
\begin{aligned}
\|F\|_{1,p}^{p}&=\mathbb{E}|F|^{p}+\mathbb{E}\|D^{X}F\|_{\mathcal{H}}^{p}.
\end{aligned}
\label{normeD1p}
\end{equation}%
\newline
We denote by $Dom (\delta^X)$ the subset of $L^{2}(\Omega, \mathcal{H})$ composed of those elements $u$ for which there exists a positive constant $c$ such that
\begin{equation}
\Big | \mathbb{E} \Big[ <D^XF, u>_{\mathcal{H}}\Big]\Big | \leqslant c\sqrt{\mathbb{E}[F^2]},\, \text{for all $F\in\mathbb{D}_{1,2}^{X}$}
\label{Domaiendelta}
\end{equation}
 For $u\in L^{2}(\Omega ;\mathcal{H}%
)$ in $Dom (\delta^X),$ $\delta ^{X}(u)$ is the element in $%
L^{2}(\Omega )$ defined by the duality relationship 
\begin{equation}
\mathbb{E}\Big[F\delta ^{X}(u)\Big]=\mathbb{E}\Big[<D_{\cdot
}^{X}F, u_{\cdot }>_{\mathcal{H}}\Big],\text{ }F\in \mathbb{D}_{1,2}^{X}.
\label{dualité}
\end{equation}

We also use the notation $\int_{0}^{T}u_{t}\delta X_{t}$ for\ $\delta
^{X}(u).$ A\ class of processes that belong to the domain of $\delta ^{X}$
is given as follows: let $\mathcal{S}^{\mid\mathcal{H}\mid}$ be the class of $%
\mathcal{H}$-valued random variables $u=\sum_{j=1}^{n}F_{j}h_{j}$ ($F_{j}\in 
\mathcal{S},$ $h_{j}\in \mid\mathcal{H}\mid).$ \\
In the same way $\mathbb{D}_{1,p}^{X}(\mid \mathcal{H\mid })$ is defined as
the completion of $\mathcal{S}^{\mid \mathcal{H\mid }}$ under the semi-norm \\  
$\parallel u\parallel _{1,p,\mid \mathcal{H\mid }}^{p}:=\mathbb{E}\parallel u\parallel_{\mid\mathcal{H}\mid}^p+\mathbb{E}\parallel D^Xu\parallel_{\mid\mathcal{H}\mid\otimes\mid\mathcal{H}\mid}^p,$ where%
\begin{equation}
\parallel D^{X}u\parallel _{\mid \mathcal{H\mid \otimes \mid H\mid }%
}^{2}=\int_{[0,T]^{4}}\mid D_{s}^{X}u_{t}\mid \mid D_{t^{\prime
}}^{X}u_{s^{\prime }}\mid \widetilde{\phi }(s, s^{\prime })\widetilde{\phi }%
(t, t^{\prime })dsdtds^{\prime }dt^{\prime }.
\label{normeDXu}
\end{equation}
The space $\mathbb{D}_{1, 2}^{X}(\mid\mathcal{H}\mid)$ is included in the domain of $\delta ^X$, and we have, for $u\in\mathbb{D}_{1, 2}^{X}(\mid\mathcal{H}\mid),$
$$\mathbb{E}\Big(\delta^X(u)^2\Big)\leqslant \mathbb{E}\parallel u\parallel_{\mid\mathcal{H}\mid}^2+\mathbb{E}\parallel D^X u\parallel_{\mid\mathcal{H}\mid\otimes\mid\mathcal{H}\mid}^2.$$

\begin{rem} Let $N_t$ be given by (\ref{N}), then we have:\\
$Var N_t=\parallel K_t^{*}\sigma\parallel_{L^{2}([0, T])} ^2= \|\sigma 1_{[0, t]} \|_{\mathcal{H}}^{2}=\displaystyle\int_{0}^{t}\displaystyle\int_{0}^{t}\phi(r, u)  \sigma_r\sigma_u du.$
 \end{rem}

\section{Proofs of the main results}
\subsection{Proof of Theorem \ref{Thmsoledsr}}
Theorem \ref{Thmsoledsr} is proven by means of an Itô formula given in \cite{DoKn}. For the convenience of the reader we state it here.\\
We have the following theorem that is proved in \cite{DoKn}.
\begin{thm}(\cite{DoKn})
Let $N_t=\int_{0}^{t}\sigma_s\delta X_s$ $(t\in[0, T]),$ and suppose that the
kernel $K$ of $X$ satisfies (H1) and (H2) and $\sigma =\{\sigma _{t}, t\in \lbrack 0,T]\}$ is a bounded function. If $F\in $ $\mathcal{C}^{1,2}([0,T]\times \mathbb{R})$
satisfies $(\ref{f})$, $\partial /\partial xF(\cdot ,N_{\cdot })\in 
\mathbb{D}_{1,2}^{X}(\mid \mathcal{H}\mathcal{\mid })$ 
and, for all $t\in \lbrack 0,T],$%
\begin{equation}
\begin{aligned}
F(t,N_{t})& =F(0,0)+\int_{0}^{t}\frac{\partial F}{\partial s}%
(s,N_{s})ds+\int_{0}^{t}\frac{\partial F}{\partial x}%
(s,N_{s})\sigma _{s}\delta X_{s}
\notag  +\frac{1}{2}\int_{0}^{t}\frac{\partial ^{2}F}{\partial
x^{2}}(s,N_{s})\frac{d}{ds}Var(N_{s})ds.
\end{aligned}
\label{ItoFormula}
\end{equation}
\label{Ito}
\end{thm}
Now, we are in position to prove Theorem \ref{Thmsoledsr}.\\

Theorem \ref{Ito} applied to $u(t, N_t)$ gives
$$du(t, N_t)=\Big(\frac{\partial u}{\partial t}(t, N_t)+\frac{1}{2}\frac{d}{dt}Var(N_t)\frac{\partial^{2} u}{\partial x^{2}}(t, N_t)\Big)dt+\sigma_t\frac{\partial u}{\partial x}(t, N_t)\delta X_t.$$
By plugging $(\ref{Edp})$ into the first term on the right side of the equation above we get $$du(t,N_t)= -f(t,N_t, u(t,N_t) ,-\sigma_t\frac{\partial u}{\partial x}(t, N_t))dt +\sigma_t\frac{\partial u}{\partial x}(t,N_t)\delta X_t.$$
 Therefore the pair $(Y, Z)$ given by $Y_t=u(t, N_t)$ and $Z_t=-\sigma_t\frac{\partial u}{\partial x}(t, N_t)$ is a 
solution of (\ref{BSDE}). 
From (\ref{f}) we conclude $\mathbb{E}\int_{t_0}^{T}Z_s^2ds<\infty$ and $\mathbb{E}\int_{t_0}^{T}Y_s^2ds<\infty.$ This proves the theorem.\,\,$\blacksquare$\\
\subsection{Proof of Theorem \ref{connexionedp}}
The proof of this theorem needs some auxiliary results.
\\
Let
\begin{equation}
p_{\mid t \mid}(x)=\frac{1}{\sqrt{2\pi \mid t \mid}}e^{-{x^{2}/2\mid t \mid}},\,\,\, t\in\mathbb{R^*}.
\label{p}
\end{equation}
Let $k$ be a continuous function, such that the following is well defined:
\begin{equation}
P_{\mid t \mid}k(x)=\displaystyle\int_{\mathbb{R}}p_{\mid t \mid}(x-y)k(y)dy.
\label{P}
\end{equation}
 A straightforward calculation shows
\begin{equation}
\frac{\partial}{\partial t}P_{\mid t \mid}k(x)=\frac{1}{2}sig(t)\frac{\partial^{2}}{\partial x^{2}}P_{\mid t \mid}k(x).
\label{relationdtdx2}
\end{equation}

The following proposition will be needed in the proof of Lemma \ref{lemmaab} below.
\begin{prop}
Let $Var(N_t)$ be increasing. Assume that $h$ is a continuous function, and suppose that there exist positive constants $c$ and $\lambda^{\prime}<(8Var(N_{T}))^{-1}$ such that $\mid h(x)\mid \leqslant ce^{\lambda^{\prime}x^2}$ for all $x\in\mathbb{R}$. Then 
\begin{eqnarray}
h(N_s)&=& P_{Var(N_s)}h(0)+\displaystyle\int_{0}^{s}\frac{\partial}{\partial x}\Big[P_{ Var(N_s)-Var(N_t)}h\Big](N_t)\sigma_t\delta X_t, 
\label{genertheoHU}
 \text{ $0<s<T.$}\nonumber
\end{eqnarray}
\label{finteg}
\end{prop}
\quad Proof. We want to apply the It\^o formula (Theorem \ref{Ito}) to $F(t,N_t)=\Big[P_{Var(N_s)-Var(N_t)}h\Big](N_t),$ $t< s.$ We begin by verifying the hypotheses of the It\^o formula:\\
\begin{align*}
\bullet \mid F(t, x)\mid &=\mid P_{ Var N_s-Var N_t }h(x)\mid=\int_{\mathbb{R}}p_{ Var N_T-Var N_t }(x-y)\mid h(y)\mid dy\\
&= \int_{\mathbb{R}}\frac{1}{\sqrt{2\pi( Var N_s-Var N_t)}}e^{-\frac{(x-y)^2}{2( Var N_s-Var N_t)}}\mid h(y)\mid dy\\
&\leqslant c\int_{\mathbb{R}}\frac{1}{\sqrt{2\pi(Var N_s-Var N_t)}}e^{-\frac{(x-y)^2}{2( Var N_s-Var N_t)}}e^{\lambda ^{\prime}y^2}dy\\
&=\frac{1}{\sqrt{2\pi(Var N_s-Var N_t) }}\int_{\mathbb{R}}e^{-\frac{(x-y)^2}{2( Var N_s-Var N_t)}}e^{\lambda^{\prime} y^2}dy\\
&=\frac{c}{\sqrt{2\pi( Var N_s-Var N_t) }}\int_{\mathbb{R}}e^{-\frac{z^2}{2(Var N_s-Var N_t)}}e^{\lambda^{\prime} (x+z)^2}dz\\
&=\frac{ce^{2\lambda^{\prime} x^2}}{\sqrt{2\pi(Var N_s-Var N_t) }}\int_{\mathbb{R}}e^{-\frac{z^2}{2(Var N_s-Var N_t)}}e^{2\lambda^{\prime} z^2}dz\\
&\leqslant\frac{ce^{2\lambda^{\prime} x^2}}{\sqrt{2\pi(Var N_s-Var N_t)}}\int_{\mathbb{R}}e^{-\frac{z^2}{2(Var N_s-Var N_t)}}e^{2\lambda^{\prime}\frac{Var N_s}{(Var N_s-Var N_t)} z^2}dz\\
&=\frac{ce^{2\lambda^{\prime} x^2}}{\sqrt{2\pi(Var N_s-Var N_t)}}\int_{\mathbb{R}}e^{-\frac{z^2}{( Var N_s-Var N_t)}(\frac{1}{2}-2\lambda^{\prime}  Var N_s)}dz\\
&\leqslant\frac{ce^{2\lambda^{\prime} x^2}}{\sqrt{2\pi( Var N_s-Var N_t)}}\sqrt{\frac{2\pi( Var N_s-Var N_t) }{1-4\lambda^{\prime} Var N_s}}=\frac{ce^{2\lambda^{\prime} x^2}}{\sqrt{1-4\lambda^{\prime} Var N_s}}\leqslant Me^{\lambda x^2},
\end{align*}
where the last inequality follows by choosing $M = \frac{c}{\sqrt{1-4\lambda Var N_T}}$ and $\lambda=2\lambda^{\prime}.$\\
The proofs for the upper bounds of $\Big | \frac{\partial}{\partial x}F(t, x)\Big |$ and $\Big | \frac{\partial^2}{\partial x^2}F(t, x)\Big |$ are similar.
\begin{align*}
\bullet \Big | \frac{\partial}{\partial x}F(t, x)\Big | 
&\leqslant\frac{2\pi c}{(2\pi(Var N_T-Var N_t) )^{\frac{3}{2}}}\int_{\mathbb{R}}\mid x-y \mid e^{-\frac{(x-y)^2}{2(Var N_s-Var N_t)}}e^{\lambda^{\prime} y^2}dy\\
&=\frac{2\pi c}{[2\pi(Var N_s-Var N_t) ]^{\frac{3}{2}}}\int_{\mathbb{R}}\mid z^{\prime} \mid e^{-\frac{(z^{\prime})^2}{2(Var N_s-Var N_t)}}e^{\lambda^{\prime} (x-z^{\prime})^2}dz^{\prime}\\
&\leqslant \frac{ce^{2\lambda^{\prime} x^2}}{2\pi(Var N_s-Var N_t)}\int_{\mathbb{R}}\frac{\mid z^{\prime} \mid }{\sqrt{ Var N_s-Var N_t} }e^{-\frac{(z^{\prime})^2}{2(Var N_s-Var N_t)}}e^{\frac{2\lambda^{\prime} Var N_s}{(Var N_s-Var N_t)} (z^{\prime}) ^2}dz^{\prime},\quad \text{($t<s$)}.
\end{align*}
Moreover, for any $\epsilon>0,$ there is a constant $K_{\epsilon}$ such that $\frac{\mid z^{\prime} \mid }{\sqrt{ Var N_s-Var N_t} }\leqslant K_{\epsilon}e^{\frac{\epsilon}{Var N_s-Var N_t}(z^{\prime})^2}.$ Then, we get
\begin{align*}
\Big | \frac{\partial}{\partial x}F(t, x)\Big | 
&\leqslant cK_{\epsilon}\frac{e^{2\lambda^{\prime} x^2}}{2\pi(Var N_s-Var N_t) }\int_{\mathbb{R}}e^{-\frac{(z^{\prime})^2}{2(Var N_s-Var N_t)}(1-2\epsilon-4\lambda^{\prime} Var N_s)}dz^{\prime}\\
&=cK_{\epsilon}\frac{e^{2\lambda^{\prime} x^2}}{\sqrt{2\pi( Var N_s-Var N_t)}}\frac{1}{\sqrt{1-2\epsilon-4\lambda^{\prime} Var N_s}}\leqslant M e^{\lambda x^2},\quad \text{(for fixed $t,$ $t<s$)}
\end{align*}
with a suitable constant $M$ and $2\lambda^{\prime}=\lambda.$
\begin{align*}
\bullet  \Big | \frac{\partial^2}{\partial x^2}F(t, x)\Big |
&\leqslant \int_{\mathbb{R}}\frac{2\pi}{[2\pi( Var N_s-Var N_t)]^{\frac{3}{2}}}e^{-\frac{(x-y)^2}{2\mid Var N_T-Var N_t\mid}}\mid h(y)\mid dy\\
&+ \int_{\mathbb{R}}\frac{2\pi\mid x-y\mid^2}{[2\pi(Var N_s-Var N_t)]^{\frac{5}{2}}}e^{-\frac{(x-y)^2}{2(Var N_s-Var N_t)}}\mid h(y)\mid dy\\
&\leqslant 2Me^{\lambda x^2}+2\pi c\int_{\mathbb{R}}\frac{(z^{\prime})^2}{[2\pi(Var N_s-Var N_t) ]^{\frac{5}{2}}}e^{-\frac{(z^{\prime})^2}{2( Var N_s-Var N_t)\mid}} e^{\lambda^{\prime} (x-z^{\prime})^2 }dz^{\prime}\\
&\leqslant 2Me^{\lambda x^2}+2\pi c\frac{e^{2\lambda^{\prime} x^2 }}{[2\pi( Var N_s-Var N_t) ]^{\frac{5}{2}}}\int_{\mathbb{R}}(z^{\prime})^2e^{-\frac{(z^{\prime})^2}{ Var N_s-Var N_t}(\frac{1-4\lambda^{\prime} Var N_s}{2})}dz^{\prime}\\
&=2Me^{\lambda x^2}+2\frac{ce^{2\lambda^{\prime} x^2 }}{(2\pi)^{\frac{3}{2}}( Var N_s-Var N_t)(1-4\lambda^{\prime}  Var N_s)^{\frac{3}{2}}}\leqslant M^{\prime}e^{\lambda x^2},\quad \text{(for fixed $t,$ $t<s$)}
\end{align*}
with a suitable constant $M^{\prime}$ and $2\lambda^{\prime}=\lambda.$\\
Moreover, for fixed $t,$ $t<s$ we have $\frac{\partial}{\partial t}P_{Var N_s-Var N_t}h(x)=\frac{1}{2}\frac{\partial^2}{\partial x^2}P_{ Var N_s-Var N_t\mid}h(x),$ and we conclude $\mid \frac{\partial}{\partial t}F(t, x)\mid \leqslant M^{"}e^{\lambda x^2}$ with a suitable constant $M^{"}.$\\
Applying the It\^o formula (Theorem \ref{Ito}) to $F(t,N_t)=\Big[P_{Var(N_s)-Var(N_t)}h\Big](N_t)$ on $[0, s-\epsilon],$ we obtain
\begin{eqnarray*}
F(s-\epsilon, N_{s-\epsilon})=F(0, 0)&+&\displaystyle\int_{0}^{s-\epsilon}\frac{\partial}{\partial x}\Big[P_{ Var(N_s)-Var(N_t)}h\Big](N_t)\sigma_t\delta X_t+\displaystyle\int_{0}^{s-\epsilon}\frac{\partial}{\partial t}\Big[P_{  Var(N_s)-Var(N_t) }h\Big](N_t)dt \\&+&\frac{1}{2}\int_{0}^{s-\epsilon}\frac{\partial^{2}}{\partial x^{2}}\Big[P_{  Var(N_s)-Var(N_t) }h\Big](N_t)\frac{d}{dt}Var(N_t)dt.
\end{eqnarray*}
Furthermore, we have
\begin{align}
\frac{\partial}{\partial t}P_{ Var(N_s)-Var(N_t)}h(x)&=\frac{\partial}{\partial \tau}P_{\tau}h(x)\mid_{\tau=Var N_s-Var N_t}\frac{\partial \tau}{\partial t}\mid_{\tau=Var N_s-Var N_t}\nonumber\\
&=\frac{1}{2}\frac{\partial^2}{\partial x^2}P_{\tau}h(x)\mid_{\tau=Var N_s-Var N_t}(-\frac{d}{dt}Var(N_t))\nonumber\\
&=-\frac{1}{2}\frac{d}{dt}Var(N_t)\frac{\partial^2}{\partial x^2}P_{Var N_s-Var N_t}h(x).
\label{rela}
\end{align}
Therefore
\begin{eqnarray*}
h(N_{s-\epsilon})-P_{Var(N_s)}h(0)&=&\displaystyle\int_{0}^{s-\epsilon}\frac{\partial}{\partial x}\big[P_{  Var(N_s)-Var(N_t) }h\big](N_t)\sigma_t\delta X_t.
\end{eqnarray*}
with $\epsilon\rightarrow 0$ we get
\begin{eqnarray}
h(N_{s})&=&P_{Var(N_s)}h(0)+\displaystyle\int_{0}^{s}\frac{\partial}{\partial x}\big[P_{  Var(N_s)-Var(N_t) }h\big](N_t)\sigma_t\delta X_t,
\label{relationendx}
\end{eqnarray}
since $h$ is a continuous function and $\frac{\partial}{\partial x}\big[P_{  Var(N_s)-Var(N_t) }h\big](N_t) \in \mathbb{D}_{1, 2}^X(\mid \mathcal{H}\mid)\subset Dom(\delta^X).$
$\blacksquare$
\begin{rems}
 If $Var N_t$ is decreasing, we have for all $0\leqslant t< s<T$ 
\begin{align*}
\frac{\partial}{\partial t}P_{\mid Var(N_s)-Var(N_t)\mid }h(x)
&=-\frac{1}{2}\frac{d}{dt}Var(N_t)\frac{\partial^2}{\partial x^2}P_{\mid Var(N_s)-Var(N_t)\mid }h(x).
\end{align*}
In this case we must add the following hypothesis in the previous proposition $\frac{Var (N_s)}{\mid Var (N_s)-Var (N_t)\mid}\geqslant 1,$ for all $0\leqslant t< s<T.$ We need this hypothesis to verify (\ref{f}) and to apply the It\^o formula.
\end{rems}
The following lemma will play an important role in this paper.
\begin{lem}
Suppose that $K_T^*$ is injective. Let $b(s, x)$ and $a(s, x)$, $s\in(t_0, T),$ $x\in\mathbb{R}$ be continuous functions, $b$ continuously differentiable with respect to $x$ and of a polynomial growth. Let $a(., N_.)\in Dom(\delta^X).$ If
\begin{equation}
\int_{t_0}^{t}b(s, N_s) ds + \int_{t_0}^{t}a(s, N_s)\delta X_s=0,\; \text{for all $t\in[t_0, T],$}
\label{equalem1}
\end{equation}
then
\begin{equation}
b(s, x)=a(s, x)=0\quad \text{for all $s\in(t_0, T),\; x\in\mathbb{R}$}.
\label{ab}
\end{equation}
\label{lemmaab}
\end{lem}
\quad Proof. First we show that
\begin{equation}
b(s, N_s)=\mathbb{E}[b(s, N_s)]+\int_{0}^{s}\Big(\int_{\mathbb{R}}\frac{\partial}{\partial x}p_{Var(N_s) -Var(N_u)}(N_u-y)b(s, y) dy\Big)\sigma_u\delta X_u.
\label{thm123}
\end{equation}
In fact, on the one hand, for all $u\in(0, s),$ we have
\begin{eqnarray*}
\displaystyle\int _{\mathbb{R}}\frac{\partial}{\partial x}p_{Var(N_s) -Var(N_u) }(N_u-y)b(s, y) dy&=&\frac{\partial}{\partial x}\displaystyle\int_{\mathbb{R}}p_{ Var(N_s) -Var(N_u) }(N_u-y)b(s, y) dy\\&=&\frac{\partial}{\partial x}P_{ Var(N_s) -Var(N_u) }b(s, N_u).
\end{eqnarray*}
Moreover, we have
\begin{eqnarray*}
\mathbb{E}\Big[b(s, N_s)\Big]
&=&\displaystyle\int_{\mathbb{R}}b(s, x)\frac{e^{-\frac{x^{2}}{2Var (N_s)}}}{\sqrt{2\pi Var(N_s)}} dx, \;\; since \quad N_s\sim \mathcal{N}(0, Var(N_s)).
\end{eqnarray*}
Therefore, 
\begin{equation*}
\mathbb{E}\Big[b(s, N_s)\Big]
=\displaystyle\int_{\mathbb{R}}b(s, x)p_{Var(N_s)}(x) dx
=\displaystyle\int_{\mathbb{R}}b(s, x)p_{Var(N_s)}(-x) dx
= P_{Var(N_s)}b(s,0).
\end{equation*}
On the other hand, we apply Proposition \ref{finteg} to $h(x)=b(s, x)$ (for fixed $s$) and
\begin{eqnarray*}
b(s, N_s)&=&P_{Var(N_s)}b(s, 0)+\int_{0}^{s}\frac{\partial}{\partial x}\Big[P_{Var(N_s) -Var(N_u) }b(s, .)\Big](N_u)\sigma_u\delta X_u.
\end{eqnarray*}
Therefore,
\begin{eqnarray*}
b(s, N_s)
&=&\mathbb{E}[b(s, N_s)]+\int_{0}^{s}\Big(\int_{\mathbb{R}}\frac{\partial}{\partial x}p_{ Var(N_s) -Var(N_u) }(N_u -y)b(s, y)dy\Big)\sigma_u\delta X_u.
\label{bs}
\end{eqnarray*}
Now, we show that
\begin{align}
&\int_{t_0}^{t}\Big[\int_{0}^{s}\Big(\int_{\mathbb{R}}\frac{\partial}{\partial x}p_{ Var(N_s) -Var(N_u) }(N_u -y)b(s, y)dy\Big)\sigma_u\delta X_u\Big]ds\nonumber\\&=\int_{t_0}^{t}\sigma_u\int_{u}^{t}\int_{\mathbb{R}}\frac{\partial}{\partial x}p_{ Var(N_s) -Var(N_u) }(N_u -y)b(s, y)dyds\delta X_u.
\label{FubStoch}
\end{align}
In fact, for $F\in L^{2}(\Omega, \mathcal{H})$
\begin{align*}
&\mathbb{E}\Big[\int_{t_0}^{t}\Big[\int_{0}^{s}\Big(\int_{\mathbb{R}}\frac{\partial}{\partial x}p_{ Var(N_s) -Var(N_u) }(N_u -y)b(s, y)dy\Big)\sigma_u\delta X_u\Big]dsF\Big]\\
&=\int_{t_0}^{t}\mathbb{E}\Big[\int_{0}^{s}\Big(\int_{\mathbb{R}}\frac{\partial}{\partial x}p_{ Var(N_s) -Var(N_u) }(N_u -y)b(s, y)dy\Big)\sigma_u\delta X_uF\Big]ds\\
&=\int_{t_0}^{t}\mathbb{E}\Big[<\int_{\mathbb{R}}\frac{\partial}{\partial x}p_{ Var(N_s) -Var(N_.) }(N_. -y)b(s, y)dy\sigma_., D_.^XF>_{\mathcal{H}}\Big]ds\\
&=\mathbb{E}\Big[<\int_{.}^{t}\Big(\int_{\mathbb{R}}\frac{\partial}{\partial x}p_{ Var(N_s) -Var(N_.) }(N_. -y)b(s, y)dy\Big)ds\sigma_., D_.^XF>_{\mathcal{H}}\Big]\\
&=\mathbb{E}\Big[\int_{t_0}^{t}\sigma_u\Big(\int_{u}^{t}\int_{\mathbb{R}}\frac{\partial}{\partial x}p_{ Var(N_s) -Var(N_u) }(N_u -y)b(s, y)dyds\Big)\delta X_uF\Big].
\end{align*}
From (\ref{thm123}) and (\ref{FubStoch}) we get
\begin{eqnarray*}
\int_{t_0}^{t}b(s, N_s)ds&=&\int_{t_0}^{t}\mathbb{E}[b(s, N_s)]ds+\int_{t_0}^{t}\Big[\int_{0}^{s}\Big(\int_{\mathbb{R}}\frac{\partial}{\partial x}p_{ Var(N_s) -Var(N_u) }(N_u -y)b(s, y)dy\Big)\sigma_u\delta X_u\Big]ds.\\
&=&\int_{t_0}^{t}\mathbb{E}[b(s, N_s)]ds+\int_{t_0}^{t}\sigma_u\int_{u}^{t}\int_{\mathbb{R}}\frac{\partial}{\partial x}p_{ Var(N_s) -Var(N_u) }(N_u -y)b(s, y)dyds\delta X_u.
\end{eqnarray*}
Then from $(\ref{equalem1})$ we get
\begin{align*}
&0=\int_{t_0}^{t}b(s, N_s) ds + \int_{t_0}^{t}a(s, N_s)\delta X_s\\
&= \int_{t_0}^{t}\mathbb{E}[b(s, N_s)]ds+\int_{t_0}^{t}a(s, N_s)\delta X_s+\int_{t_0}^{t}\sigma_u\Big[\int_{u}^{t}\Big(\int_{\mathbb{R}}\frac{\partial}{\partial x}p_{ Var(N_s) -Var(N_u) }(N_u -y)b(s, y)dy\Big)ds\Big]\delta X_u\\
&=\int_{t_0}^{t}\Big[a(u, N_u)+\sigma_u\int_{u}^{t}\int_{\mathbb{R}}\frac{\partial}{\partial x}p_{ Var(N_s) -Var(N_u) }(N_u -y)b(s, y)dy ds\Big]\delta X_u+\int_{t_0}^{t}\mathbb{E}[b(s, N_s)]ds.
\end{align*}
Thus
\begin{eqnarray}
\int_{t_0}^{t}\mathbb{E}[b(s, N_s]ds&=&0,
\label{eq6}
\end{eqnarray}
and
\begin{eqnarray}
\int_{t_0}^{t}\Big[a(u, N_u)+\sigma_u\int_{u}^{t}\int_{\mathbb{R}}\frac{\partial}{\partial x}p_{ Var(N_s) -Var(N_u) }(N_u -y)b(s, y)dy ds\Big]\delta X_u&=&0.
\label{eq7}
\end{eqnarray}
Let
\begin{eqnarray*}
Z(u, N_u)&:=&a(u,N_u)+\sigma_u\int_{u}^{t}\int_{\mathbb{R}}\frac{\partial}{\partial x}p_{ Var(N_s) -Var(N_u) }(N_u-y)b(s, y)dy ds
\end{eqnarray*}
(\ref{eq7}) implies that $Var(\int_{t_0}^{t} Z(u, N_u)\delta X_u)=0.$ Therefore $K_t^*Z(., N_.)_u=0$ $u-$a.e. in $[t_0, t].$\\
By injectivity of $K_t^*,$ we get
\begin{eqnarray*}
Z(u, N_u)=0,\, \text{for all $u\in[t_0, t]$.}
\end{eqnarray*}
Therefore
\begin{eqnarray*}
a(u, z)+\sigma_u\int_{u}^{t}\int_{\mathbb{R}}\frac{\partial}{\partial x}p_{ Var(N_s) -Var(N_u) }(z -y)b(s, y)dy ds=0,
\end{eqnarray*}
for all $z\in\mathbb{R}$ (\cite{HuPeng}).
Now by differentiating with respect to $t$, we get
\begin{eqnarray*}
\int_{\mathbb{R}}\frac{\partial}{\partial x}p_{ Var(N_t) -Var(N_u)}(z -y)b(t, y)dy=0.
\end{eqnarray*}
for all $t>u$ and $z\in\mathbb{R}.$\\
An integration by parts formula yields 
\begin{eqnarray*}
\int_{\mathbb{R}}p_{Var(N_t) -Var(N_u)}(z -y)\frac{\partial}{\partial y}b(t, y)dy=0.
\end{eqnarray*}
Let $u\rightarrow t$ then we see that
$\frac{\partial}{\partial y}b(t, y)=0$ for all $t\in(t_0, T-\epsilon)$ and $y\in\mathbb{R}.$ This means that there is a $b_1(t)$ such that $b(t, y)=b_{1}(t).$
Now from $(\ref{eq6})$ we have
\begin{eqnarray*}
\int_{t_0}^{t}b_{1}(s) ds=0, \; t_0\leqslant t\leqslant T-\epsilon.
\end{eqnarray*}
This implies that $b_1(t)=0$ for all $t\in(t_0, T)$ and accordingly $b(t, y)=0$ for all $t\in(t_0, T)$ and  $y\in\mathbb{R}$. Thus, $a(t, y)=0$.\,\,$\blacksquare$\\
\begin{rems}
In the proof of Lemma \ref{lemmaab} we have $$Var\int_{t_{0}}^{t}Z(u, N_{u})\delta
X_{u}=Var\int_{t_{0}}^{t}K_{t}^{\ast }Z(\cdot, N_{\cdot })_{u}\delta W_{u}=0$$
which implies that $K_{t}^{\ast }Z(\cdot ,N_{\cdot })_{u}=\int_{u}^{t}\frac{%
\partial }{\partial v}K(v, u)Z(v, N_{v})dv=0$ $\ (u, \omega )-a.e.$ on $%
[t_{0}, t]\times \Omega ,$ for all $t\in(t_0, T].$ We would like to conclude that $Z(v, N_{v})=0$ $\ (v, \omega )-a.e.$ on $%
[t_{0}, t]\times \Omega .$\\
An evident hypothesis is the injectivity of $%
K_{t}^{\ast }Z(\cdot, N_{\cdot })_{u}$ as a function of $u\in \lbrack
t_{0}, t].$ \\
Let us look for a sufficient condition for injectivity:

\[
0=\int_{t_{0}}^{t}K_{t}^{\ast }Z(\cdot ,N_{\cdot
})_{u}du=\int_{t_{0}}^{t}Z(s, N_{s}, \omega
)\int_{t_{0}}^{s}\frac{\partial }{\partial s}K(s, u)duds.
\]%
Suppose that $K_{t}^{\ast }Z(\cdot, N_{\cdot })_{u}=0$ $\ u-a.e.$ on $%
[t_{0}, t]$ for all $t\in (t_{0}, T]$ and $Z>0$ on $%
(a,b)\subset \lbrack t_{0}, T].$ Then
\begin{eqnarray*}
0 &=&\int_{t_{0}}^{b}K_{b}^{\ast }Z(\cdot, N_{\cdot
})_{u}du-\int_{t_{0}}^{a}K_{a}^{\ast }Z(\cdot, N_{\cdot })_{u}du \\
&=&\int_{t_{0}}^{b}Z(s, N_{s}, \omega )\int_{t_{0}}^{s}\frac{\partial }{%
\partial s}K(s, u)duds-\int_{t_{0}}^{a}Z(s, N_{s}, \omega )\int_{t_{0}}^{s}%
\frac{\partial }{\partial s}K(s, u)duds \\
&=&\int_{a}^{b}Z(s, N_{s}, \omega )\int_{t_{0}}^{s}\frac{\partial }{\partial s}%
K(s, u)duds
\end{eqnarray*}
Let $\widetilde{K}_{t_{0}}(s)=\int_{t_{0}}^{s}\frac{\partial }{\partial s}%
K(s, u)du.$ If we suppose that $\widetilde{K}_{t_{0}}(s)\neq 0$ on $[t_{0}, T]
$ (and therefore does not change sign on $[t_{0}, T]),$ we obtain a contradiction. Therefore, a sufficient condition for $K_T^*$ to be injective is $\widetilde{K}_{t_{0}}(s)\neq 0$ on $[t_0, T].$ This last hypothesis is satisfied in particular if $\frac{\partial }{\partial s}%
K(s, u)>0$ (or $<0)$ for $u\in (t_{0}, s)$ for all $s\in(t_{0}, T).$ 
\label{reminjectivite}
\end{rems}

We are now ready to prove Theorem \ref{connexionedp}. By Theorem \ref{Ito} we get 
\begin{eqnarray*}
du(t, N_t)
&=& \Big[\frac{\partial}{\partial t}u(t, N_t)+\frac{1}{2}\frac{d}{dt}Var(N_t) \frac{\partial^2}{\partial x^2}u(t, N_t) \Big ] dt+\sigma_t \frac{\partial}{\partial x}u(t, N_t)\delta X_t.
\end{eqnarray*}
Moreover, we can write
\begin{eqnarray*}
u(t, N_t)-g(N_T)= -\int_{t}^{T}\Big[\frac{\partial}{\partial s}u(s, N_s)+\frac{1}{2}\frac{d}{ds}Var(N_s)\frac{\partial^2}{\partial x^2}u(s, N_s)\Big]ds-\int_{t}^{T}\sigma_s \frac{\partial}{\partial x}u(s, N_s) \delta X_s.
\end{eqnarray*}
Using (\ref{BSDE}), we get 
$$
\begin{array}{ll}
&\displaystyle\int_{t}^{T}f(s, N_s, u(s, N_s), v(s, N_s))ds +\int_{t}^{T}v(s, N_s)\delta X_s\\&=
-\displaystyle\int_{t}^{T}\Big[\frac{\partial}{\partial s}u(s, N_s)+\frac{1}{2}\frac{d}{ds}Var(N_s) \frac{\partial^2}{\partial x^2}u(s, N_s)\Big]ds-\int_{t}^{T}\sigma_s \frac{\partial}{\partial x}u(s, N_s) \delta X_s.
\end{array}
$$
We evaluate for $t = t_0$ and we make the subtracting with the above equation.
We obtain for all $t_0\leqslant t\leqslant T-\epsilon$
$$
\begin{array}{l}
\displaystyle\int_{t_0}^{t}f(s, N_s, u(s, N_s),  v(s, N_s))ds +\displaystyle\int_{t_0}^{t}v(s, N_s)\delta X_s\\
=-\displaystyle\int_{t_0}^{t}\Big[\frac{\partial}{\partial s}u(s, N_s)+\frac{1}{2}\frac{d}{ds}Var(N_s)\frac{\partial^2}{\partial x^2}u(s, N_s)\Big]ds-\displaystyle\int_{t_0}^{t}\sigma_s \frac{\partial}{\partial x}u(s, N_s) \delta X_s.
\end{array}
$$
Using (\ref{Edp}), we obtain
\begin{align*}
&\displaystyle\int_{t_0}^{t}\Big[f(s, N_s, u(s, N_s),  v(s, N_s))-f(s, N_s, u(s, N_s),  -\sigma_s \frac{\partial}{\partial x}u(s, N_s))\Big]ds \\&+\displaystyle\int_{t_0}^{t}\Big[v(s, N_s)+\sigma_s \frac{\partial}{\partial x}u(s, N_s)\Big] \delta X_s=0.
\end{align*}
Since $u,$ $v$ satisfy (\ref{f}) and $f\in\mathcal{C}^{0, 1}([0, T ] \times  \mathbb{R}^3)$ is of polynomial growth, we apply now Lemma \ref{lemmaab} and obtain
\begin{eqnarray*}
v(t, x)=-\sigma_{t} \frac{\partial}{\partial x}u(t, x), \; \forall t\in(t_0, T), x\in\mathbb{R}.\,\,\text{$\blacksquare$}
\end{eqnarray*}
\subsection{Proof of Theorem \ref{Comparisontheo} and Theorem \ref{continuityY}}

We start with proving this preliminary result:
\begin{prop}(Transfer formula)
Let $X$ be given by (\ref{X}). Let $D^W$ be the Malliavin derivative with respect to the Brownian motion. Then, $K_T^*D_.^X=D_.^W$ on $\mathbb{D}_{1, 2}^{X}.$
\label{Transferformula}
\end{prop}
Proof. We have $F=f(X(\varphi_1), ..., X(\varphi_n))=f(W(K_T^*\varphi_1), ..., W(K_T^*\varphi_n)),$ $n\geqslant 1.$\\
Recall that $X(\varphi_i)=\int_{0}^{T}\varphi_i(s) \delta X_s= \int_{0}^{T}(K_T^*\varphi_i)_s\delta W_s=W(K_T^*\varphi_i), i=1,..., n, ,  \varphi_1,..., \varphi_n\in\mathcal{H}.$\\
Furthermore
\begin{align*}
(K_T^{*}D_.^X F)_t&=\displaystyle\int_ {t}^{T}D_u^X F\frac{\partial K}{\partial u}(u, t) du\\&=\sum\limits_{i=1}^{n}\frac{\partial f}{\partial x_i}(X(\varphi_1), ..., X(\varphi_n))\int_{t}^{T}\varphi_i(u)\frac{\partial K}{\partial u}(u, t) du\\
&=\sum\limits_{i=1}^{n}\frac{\partial f}{\partial x_i}(X(\varphi_1), ..., X(\varphi_n))(K_T^*\varphi _i)_t=D_t^W F, \text{$t\in[0, T]$}.
\end{align*}
We extend the equality $(K_T^{*}D_.^X F)_t =D_t^W F$ to the closure of the linear combinations of $\mathcal{S}$ by means of the norm
\begin{equation*}
\begin{aligned}
\|F\|_{1,2}^{2}&=\mathbb{E}|F|^{2}+\mathbb{E}\|D^{X}F\|_{ \mathcal{H}}^{2}.
\end{aligned}
\end{equation*}%
We have $F_1, F_2\in \mathcal{S}, a_1, a_2\in\mathbb{R}:$ $(K_{T}^{*}D_{.}^{X}(a_1 F_1+a_2F_2))_t=D_t^{W}(a_1 F_1+a_2F_2)$ because $K^*, D^X$ and $D^W$ are linear operators.\\
Let $(F_n)_{n\in \mathbb{ N}}\subset span(\mathcal{S})$ such that $(F_n)_{n\in \mathbb{ N}}$ converges in norm $\|.\|_{1,2}$ to $F.$ Then, $F\in \mathbb{D}_{1,2}^{X}.$ Thus, $\| D^X F- D^X F_n \|_{ \mathcal{H}}\rightarrow 0,$ with $n \longmapsto \infty.$ Thus, $\| K_T^*D^X F- K_T^*D^X F_n \|_{L^{2}(\Omega, L^{2}(0, T))}\rightarrow 0,$ with $n \longmapsto \infty$ by isometry.\\
Since $K_T^*D^XF_n=D^WF_n,$ $(D^W F_n)_{n\in\mathbb{N}}$ is a convergent sequence in $L^{2}(\Omega, L^{2}(0, T)).$ By proposition 1.2.1 (\cite{Nualart}) the limit in $L^{2}(\Omega, L^{2}(0, T))$ is $D^WF.$ Therefore, $K_T^*D_.^X=D_.^W.\,\,\blacksquare $ \\
\subsubsection{Proof of Theorem \ref{Comparisontheo}}
Let $u^i(t, x)$ be the solution to (\ref{Edp}) with $f$ replaced by $f^i, g$ replaced by $g^i.$ Then the solution to (\ref{edsrmulti}) is given by $Y^i_t=u^i(t, N_t).$ It suffices to prove $u^1(t, x)\geqslant u^2(t, x).$\\
Denote $\rho_t=\sqrt{\frac{d}{dt}Var(N_t)},$ and $\zeta_t=\int_{0}^{t}\rho_s\delta W_s,$ where $W$ is standard Brownian motion.\\
Applying It\^o's formula with respect to $W,$ we have
\begin{align*}
du^i(t, \zeta_t)&=\rho_t\frac{\partial u^i}{\partial x}(t, \zeta_t)dW_t+\frac{\partial u^i}{\partial t}(t, \zeta_t)dt+\frac{1}{2}\frac{\partial^2 u^i}{\partial x^2}(t, \zeta_t)\frac{d}{dt}Var(N_t)dt\\
&=-f^i(t, \zeta_t, u^i(t, \zeta_t), -\sigma_t\frac{\partial u^i}{\partial x}(t, \zeta_t))dt+\rho_t\frac{\partial u^i}{\partial x}(t, \zeta_t)dW_t.
\end{align*}
Thus, $(\tilde{Y}^i_t, \tilde{Z}^i_t)=(u^i(t, \zeta_t), \rho_t\frac{\partial u^i}{\partial x}(t, \zeta_t))$ is a solution to the following BSDE
\begin{equation*}
\begin{cases}
 d\tilde{Y}^i_t=-f^i(t, \zeta_t, \tilde{Y}^i_t, -\sigma_t\rho_t^{-1}\tilde{Z}^i_t)dt+\tilde{Z}^i_tdW_t, &\text{$0<t_{0}\leqslant t\leqslant T$}\\
\tilde{Y}^i_T=g^i(\zeta_T)\\
\end{cases}
\end{equation*}
By the classical comparison theorem (\cite{ZhangI}), $\tilde{Y}^1_t\geqslant \tilde{Y}^2_t$ almost surely. Thus, $u^1(t, \zeta_t)\geqslant u^2(t, \zeta_t).$ Since $\zeta_t$ is a gaussian random variable with positive variance, from lemma 3.7 (\cite{HuPeng}), we conclude $u^1(t, x)\geqslant u^2(t, x).$\,\,$\blacksquare$

\subsubsection{Proof of Theorem \ref{continuityY}}
 The proof of this theorem is based on an application of theorem 5.1 in \cite{D.Nualart}. In fact, we have
\begin{align*}
\parallel D_.^XY_t\parallel_{\mathcal{H}}&=\parallel (K_T^*D^XY)_t\parallel_{L^2[0, T]}\\
&=\parallel D_.^WY_t\parallel_{L^{2}[0, T]},\quad \text{(by Proposition \ref{Transferformula}).}
\end{align*}
Moreover
$$D_s^WY_t=D_s^Wu(t, N_t)=\frac{\partial}{\partial x}u(t, N_t)(K_t^*\sigma)_s.$$
Therefore
\begin{align*}
\parallel D_.^XY_t\parallel_{\mathcal{H}}&=\parallel \frac{\partial}{\partial x}u(t, N_t)(K_t^*\sigma)_.\parallel_{L^2[0, T]}\\
&=\int_{0}^{T}\Big(\frac{\partial}{\partial x}u(t, N_t)\Big)^2(K_t^*\sigma)_s^2ds\\
&=\Big(\frac{\partial}{\partial x}u(t, N_t)\Big)^2Var N_t.
\end{align*}
Under the hypothesis $Var N_t>0$ and $\frac{\partial}{\partial x}u(t, N_t)\neq 0,$ we apply theorem 5.1 in \cite{D.Nualart}, we get \\$\parallel D_.^XY_t\parallel_{\mathcal{H}}>0.$\,\,$\blacksquare$

\section*{Acknowledgments}
I am grateful to M. Marco Dozzi for stimulating discussions on this topic. I thank the program Hubert Curien "Utique" of the 'French Ministry of Foreign Affairs' and the 'Tunisian Ministry of Education and Research' for the financial support. \newline

\bibliographystyle{plain}
\bibliography{KnaniDozzibiblio3}

\end{document}